\documentclass[12pt]{article}
\usepackage{latexsym}
\usepackage{amsfonts}
\def\ZZ{\mathbb Z}
\def\RR{\mathbb R}
\def\CC{\mathbb C}

\def\CP{\mathbb C \mathbb P}
\def\BW{\mathcal B \mathcal W}
\def\eea{\end{eqnarray*}}

\newtheorem{defn}{Definition}
\newtheorem{thm}{Theorem}
\newtheorem{prop}[thm]{Proposition}
\newtheorem{cor}[thm]{Corollary}

\newtheorem{question}[thm]{Question}
\newenvironment{proof}{\medskip \noindent
{\bf Proof.}}{\hfill \rule{.5em}{1em}
\\}

\newenvironment{xpl}{\mbox{ }\\ {\bf  Example}\mbox{ }}{
\hfill $\diamondsuit$\mbox{}\bigskip}

\begin{document}
\sloppy
\title{Spin Manifolds, Einstein Metrics, \\
and Differential Topology}

\author{Masashi Ishida and Claude LeBrun\thanks{Supported 
in part by  NSF grant DMS-0072591.}  
  }

\date{July 14, 2001\\
Revised October 28, 2001}
\maketitle

\begin{abstract}
We show that there exist smooth, simply connected,
four-dimensional 
{\em spin}
manifolds which do not admit Einstein metrics, but nevertheless
satisfy the strict Hitchin-Thorpe inequality.
Our construction makes use of the  Bauer/Furuta 
cohomotopy refinement of the Seiberg-Witten invariant
\cite{baufu,bauer2}, in conjunction with  curvature 
estimates previously proved by the second author \cite{lric}. 
These methods also  easily allow one to  construct 
examples
of topological $4$-manifolds which admit an Einstein
metric for one smooth structure, but which have
infinitely many other smooth structures for which no
Einstein metric can exist. 
\end{abstract}

\section{Introduction}\label{sec:intro}

A smooth Riemannian metric $g$ is said to be {\em  Einstein} 
if its Ricci curvature, considered as a function 
on the unit tangent bundle,  is constant. 
As recently as the  1960's,  it could have seemed 
reasonable to
hope that every smooth compact simply connected manifold might admit
such a metric;  and, indeed,   it was apparently with this goal in mind  that 
Yamabe \cite{yamabe}  carried out his trail-blazing work 
on the total  scalar curvature. 
In the late 1960's, however, Thorpe \cite{tho}  
observed, as a parenthetical remark,  that a compact oriented $4$-dimensional 
Einstein manifold $(M^{4},g)$ must satisfy the inequality 
\begin{equation}
2\chi (M) \geq 3|\tau (M)|,	
	\label{eq:ht}
\end{equation}
 where $\chi$ denotes the 
Euler characteristic and $\tau$ denotes the signature. 
This inequality was rediscovered several years later by 
Hitchin \cite{hit}, who went on to give
examples of simply connected $4$-dimensional manifolds 
which violate inequality (\ref{eq:ht}), and thus do
not admit Einstein metrics. In light of this,   (\ref{eq:ht}) 
 has   come to be known 
as the {\em Hitchin-Thorpe inequality} \cite{bes}.

 Perhaps the deepest results in Hitchin's paper 
 concern 
 the boundary case of inequality (\ref{eq:ht}); 
 for it was  shown in  \cite{hit} 
that  a simply connected $4$-dimensional
Einstein manifold with $2\chi = -3\tau$ must be 
Ricci-flat, K\"ahler, and 
diffeomorphic to the Kummer surface $K3$.
In particular, this implies that most 
simply connected $4$-manifolds satisfying
$2\chi = 3 |\tau|$ do not admit Einstein metrics; 
and this  applies not only to  examples
of the `wrong' homotopy type --- 
e.g.   $\CP_{2}\# 9 \overline{\CP}_{2}$ ----  
 but also \cite{poonmail} to an infinite
class  of manifolds \cite{frmo,fsicm}
 now known  to be 
homeomorphic but not diffeomorphic to $K3$; cf. \cite{kot}. 
Since  Yau \cite{yau}  showed that 
 $K3$ actually admits Einstein metrics of the 
 kind described in  Hitchin's result,
 this shows, in particular,  that the existence of 
 Einstein metrics on a $4$-manifold is a matter of   
 diffeotype, and not just of homeotype. 
 
 By contrast, for simply connected $4$-manifolds satisfying the 
 so-called {\em strict Hitchin-Thorpe inequality} 
 \begin{equation}
2\chi (M) > 3|\tau (M)|,	
	\label{eq:strict}
\end{equation}
 obstructions to the existence of Einstein metrics
remained unknown until 1995, when the second author 
\cite{lno} used scalar curvature estimates derived
from the Seiberg-Witten equations to 
display an infinite class of $4$-manifolds which do not 
admit Einstein metrics,  but nonetheless 
satisfy (\ref{eq:strict}). However, the examples given there 
were  non-minimal complex surfaces of general type, 
 and thus, in particular,  were all {\em non-spin}.
Even though these results have been subsequently improved
\cite{lric} via the introduction of Weyl-curvature
terms into the estimates, the non-spin
character of the examples 
seemed to reflect an essential feature 
of the 
construction. This made it 
natural to wonder about the following:

\begin{question}\label{pair}
	Are there simply connected $4$-dimensional {\bf spin} 
	manifolds which do not admit Einstein metrics,
	but which nonetheless
	 satisfy the strict Hitchin-Thorpe inequality 
	(\ref{eq:strict})? 
	\end{question}

In Theorem \ref{spin} below, we will see  that the 
answer is affirmative. Our construction depends  
crucially on a recent breakthrough in Seiberg-Witten theory, due to 
Furuta and Bauer \cite{baufu,bauer2}, used 
in conjunction with 
recent curvature estimates due to the second author \cite{lric}.

	Another curious defect of the examples constructed  
	in \cite{lno,lebweyl,lric} is that, when oriented
	so that the signature is negative, all have the 
	the property that $b_{+}$ is odd. This merely reflects
	the fact that these examples are all almost-complex, and 
	the integrality of the Todd genus implies that 
	$b_{+}$ must be odd for any simply connected almost-complex
	$4$-manifold. However, the 
	Einstein condition is independent of orientation, and 
	these examples can certainly have $b_{-}$ of either parity. 
     Might it not be reasonable to hope that this geographical 
     curiosity was  merely  
	 an artifact  of the construction? 

\begin{question}\label{geo}
		Are there simply connected $4$-manifolds, 
		with 
		$\tau < 0$ and $b_{+}$ {\bf even}, 
		which do not admit Einstein metrics,
      but which nonetheless
	 satisfy the strict Hitchin-Thorpe inequality 
	(\ref{eq:strict})? 	
	\end{question}

In Theorem \ref{molti} below, we will see that the present methods 
allow one to construct such examples in considerable abundance. 

 Finally, we return to the heavy dependence of these
 questions on the choice of smooth structure.
 As already noted, while $K3$ admits a smooth structure for
 which Einstein metrics exist, it also admits infinitely many
 smooth structures for which Einstein metrics {\em do not}
 exist. Is this typical? 
 
 \begin{question}\label{change}
	 Let $M$ be a  compact simply connected topological $4$-manifold which 
	 admits at least one smooth structure. Are there always  
	 {\bf infinitely many} distinct smooth structures on  $M$ for which 
	 no compatible Einstein metric exists? 
	 \end{question}
	 
	 While we are certainly not able to answer this question
	 in  full generality, we are, in Theorems 
	 \ref{spin}, \ref{molti} and \ref{lots},
	  at least able
	 to give an affirmative answer for infinitely many
	 homeotypes satisfying (\ref{eq:strict}). Interestingly, 
	 many of these examples also have the property that,
	 like
	 $K3$, they also admit  smooth structures for which  compatible Einstein
	 metrics {\em do} exist. For some earlier results and ruminations 
          in this direction, see \cite{kot}. 

The infinite-fundamental-group analogue of 	 
Question \ref{change}  would still appear to be 
a  topic worthy of investigation. By contrast,
however,  the infinite-fundamental-group analogues 
of  Questions \ref{pair}
and \ref{geo} are susceptible to  the homotopy-theoretic 
volume estimates pioneered by Gromov, and  affirmative answers 
\cite{grom,samba} to
the corresponding  questions   have therefore
already been known for some time.

\section{Monopole Classes and Curvature}\label{sec:curv}

If $M$ is a smooth  oriented $4$-manifold, we can always find
Hermitian line bundles $L\to M$ such that $c_{1}(L)\equiv w_{2}(TM)
\bmod 2$. For any such $L$, and any Riemannian metric $g$ on $M$,
one can then find rank-$2$ Hermitian vector bundles ${\mathbb V}_\pm$
which formally satisfy
$${\mathbb V}_\pm= {\mathbb S}_{\pm}\otimes L^{1/2},$$
where ${\mathbb S}_{\pm}$ are the locally defined left- and
right-handed spinor bundles of $(M,g)$. 
Such a choice of   ${\mathbb V}_{\pm}$, up to isomorphism, 
is called a 
spin$^{c}$ structure $\frak c$ on $M$, and is determined, modulo
the $2$-torsion subgroup of $H_{1}(M,\ZZ)$, by the first Chern class
$c_{1}(L)= c_{1}({\mathbb V}_{\pm})
\in H^{2}(M,\ZZ )$ of the spin$^{c}$ structure. 
Because ${\mathbb S}_{+}$ is a  quaternionic line
bundle, there is a canonical anti-linear isomorphism
\begin{eqnarray*}
	{\mathbb V}_{+} & \longrightarrow & {\mathbb V}_{+}\otimes L^{*}  \\
	\Phi & \mapsto  & \bar{\Phi}
\end{eqnarray*}
called `complex conjugation.'  More importantly, every unitary 
connection $A$ on $L$ induces a Dirac operator
$$D_{A}: \Gamma ({\mathbb V}_{+})\to \Gamma ({\mathbb V}_{-}).$$
If $A$ is such a connection, and if $\Phi$ is a section of
${\mathbb V}_{+}$, the pair $(\Phi , A)$ is said to satisfy the
 {\em Seiberg-Witten equations} \cite{witten} if  
\begin{eqnarray} D_{A}\Phi &=&0\label{drc}\\
 F_{A}^+&=&-\frac{1}{2}\Phi \odot \bar{\Phi}.\label{sd}\end{eqnarray}
Here 
$F_{A}^+$  is the self-dual part of
the curvature of $A$,
and we have identified $\Lambda^{+}\otimes \CC$ with 
$[\odot^{2}{\mathbb V}_{+}]\otimes L^{*} = \odot^{2}{\mathbb S}_{+}$
  in the canonical manner. 
 
For the $4$-manifolds of primary interest here, there turn out to be 
certain spin$^{c}$ structures  for which 
there exists a solution of the Seiberg-Witten equations 
for each metric $g$. This situation is neatly codified by the 
following 
 terminology \cite{K}:
 
\begin{defn}
Let $M$ be a smooth compact oriented $4$-manifold
with $b_{+}\geq 2$. An element $a\in  H^{2}(M,\ZZ )/
\mbox{\rm torsion}$ will be called a {\bf monopole
class} of $M$ iff there exists a spin$^{c}$ structure
$\mathfrak c$ 
on $M$ with first Chern class 
$$c_{1}(L)\equiv a ~~~\bmod \mbox{\rm torsion}$$ which has the property
that  the corresponding  Seiberg-Witten 
equations (\ref{drc}--\ref{sd})
have a solution for every Riemannian  metric $g$ on $M$. 
\end{defn}

Because  the Seiberg-Witten
equations imply a uniform   {\em a priori} 
 bound  on
$|F_{A}^{+}|$ for any given metric on $g$ on $M$,
it follows \cite{witten} that  the set $\mathfrak C$ of
monopole classes of any $4$-manifold $M$ is necessarily finite.
Also note that $a\in {\mathfrak C}$ 
$\Longleftrightarrow$ 
$(-a) \in {\mathfrak C}$, since complex
conjugation sends solutions of the Seiberg-Witten 
equations to solutions of the Seiberg-Witten 
equations. Finally, notice that, for any $a,b \in  {\mathfrak C}$,
$a-b$   is  automatically
divisible by $2$ in the lattice $H^{2}(M,\ZZ )/
\mbox{\rm torsion}$, since the corresponding first Chern classes
are  both   $\equiv w_{2} \bmod 2$. 
With these observations in mind, we will 
now introduce  a crude but effective diffeomorphism invariant
of $4$-manifolds.

\begin{defn}\label{band}
	Let $M$ be a smooth compact oriented $4$-manifold 
	with $b_{+}(M)\geq 2$. Let ${\mathfrak C}\subset H^{2}(M, \ZZ )
	/\mbox{\rm torsion}$ 
	be the set of
monopole classes of $M$.  
If ${\mathfrak C}$ contains a non-zero
element (and hence at least two elements), we define
the {\bf bandwidth} of $M$ to be 
$$\BW (M )  ~=  
\max ~ \left\{ n\in \ZZ^{+}~|~ \exists ~a,b \in {\mathfrak C}, ~a\neq 
b, 
~{\rm s.t.}~ 
2n|(a-b) ~ \right\} .$$
If, on the other hand,  ${\mathfrak C}\subset \{ 0\},$
we define the bandwidth $\BW (M )$
 to be $0$. 
	\end{defn}

As was first pointed out by Witten \cite{witten},
the existence of a monopole class implies an  
{\em a priori} lower bound on the $L^{2}$ norm of 
the scalar curvature of Riemannian metrics. More
recently, the second author then discovered \cite{lebweyl,lric} that 
(\ref{drc}--\ref{sd}) also imply a family of
analogous estimates involving the 
self-dual Weyl curvature. For our
present purposes, the most useful of these estimates is the 
following  \cite{lric}:

 \begin{prop}[LeBrun]\label{est}
Let $M$ be a   smooth compact  
oriented  $4$-manifold with monopole class 
$a\in H^{2}(M,\ZZ)/\mbox{\rm torsion}\subset H^{2}(M,\RR )$. 
Let $g$ be any Riemannian metric on $M$,
and let $a^{+}\in H^{2}(M,\RR )$ be the 
self-dual part of $a$ with respect to the 
`Hodge' decomposition 
$$H^{2}(M,\RR )= {\mathcal H}^{+}_{g}\oplus {\mathcal H}^{-}_{g}$$
of second de Rham cohomology, identified with the
space of $g$-harmonic $2$-forms, into eigenspaces of the 
$\star$ operator. 
Then the scalar curvature $s$ and self-dual Weyl curvature
$W_{+}$ of $g$ satisfy 
$$
\frac{1}{4\pi^{2}}\int_{M}\left( \frac{s^{2}}{24} + 
2|W_{+}|^{2}\right) d\mu \geq \frac{2}{3} 
(a^{+})^{2} ,
$$
where $d\mu$ denotes the Riemannian volume form of 
$g$, and where the point-wise norms are calculated with
respect to $g$. Moreover, if 
$a^{+}\neq 0$, and if $a$ is not the first Chern class
of a symplectic structure on $M$, then the inequality 
is necessarily strict. 
\end{prop}

\section{The Bauer-Furuta Invariant}\label{sec:bauer}

The usual Seiberg-Witten invariant \cite{taubes3,ozsz} 
of a smooth compact oriented 
 $4$-manifold  with $b_{+}\geq 2$,
 equipped with  some   
fixed spin$^{c}$ structure $\mathfrak c$, 
 is obtained 
considering the moduli space of solutions 
of (\ref{drc}--\ref{sd}), for a generic metric $g$, 
as  a 
homology class in the configuration space 
$${\cal B}=\left( 
[\Gamma ({\mathbb V}_{+})-0]\times 
\{ \mbox{smooth unitary connections on }L ~\}\right) /{\mathcal G} ,$$
where ${\mathcal G}= \{ u: M\stackrel{C^{\infty}}{\longrightarrow} S^{1}\}$. 
However, one may instead think of this moduli space
as defining a framed bordism class in the 
space ${\cal D} \subset {\cal B}$
of solutions of (\ref{drc}). Pursuing this idea, 
Furuta and 
Bauer \cite{baufu,bauer2}
were independently able to define a refinement of the Seiberg-Witten
invariant which takes values in a cohomotopy group;
for example, if $b_{1}(M)=0$, the invariant takes 
values in $\pi^{b_{+}(M)-1}(\CP_{d-1}),$
where $d=[c_{1}^{2}(L)-\tau(M)]/8$. 
If this Bauer-Furuta stable cohomotopy  
invariant is non-zero, $c_{1}(L)$ 
is a monopole class. 
But, remarkably, this 
 invariant  is not killed off  by the sort of connect sum 
operation 
 that  would eliminate  the usual Seiberg-Witten 
 invariant \cite{bauer2}:

\begin{thm}[Bauer]
	Let $X$, $Y$, and $Z$ be compact oriented
	$4$-manifolds with $b_{1}=0$ and $b_{+}\equiv 3\bmod 4$.
	Suppose, that ${\mathfrak c}_{X}$, ${\mathfrak c}_{Y}$,
	and ${\mathfrak c}_{Z}$ are spin$^{c}$ structures 
	of almost-complex type on 
	$X$, $Y$, and $Z$ for which the mod-$2$ Seiberg-Witten 
	invariant is non-zero. Then the induced
	spin$^{c}$ structures ${\mathfrak c}_{X\# Y}$ and 
	${\mathfrak c}_{X\# Y\# Z}$ on $X\# Y$ and 
	$X\# Y \# Z$ have non-zero Bauer-Furuta stable cohomotopy 
	invariant.
	\end{thm}
	
	Now, recall that  a celebrated result of Taubes \cite{taubes} 
	asserts that the mod-$2$  Seiberg-Witten invariant 
	is non-zero for the 
	 canonical
	spin$^{c}$ structure of any  symplectic $4$-manifold.
	  Since the non-vanishing of the 
	Bauer-Furuta invariant forces the existence of 
	 solutions of the Seiberg-Witten equations, we immediately
	 obtain the following consequence: 
	
	\begin{cor}\label{voila}
			Let $X$, $Y$, and $Z$ be compact oriented
			simply connected symplectic 
	$4$-manifolds with  $b_{+}\equiv 3\bmod 4$.
	Then, with respect to the canonical isomorphisms
	\begin{eqnarray*}
		H^{2}(X\# Y, \ZZ) & = & H^{2}(X, \ZZ)\oplus H^{2}(Y, \ZZ)  \\
			H^{2}(X\# Y\# Z, \ZZ) & = &  H^{2}(X, \ZZ)\oplus H^{2}(Y, \ZZ) 
			\oplus H^{2}(Z, \ZZ) , 
	\end{eqnarray*}
	the cohomology classes 
	\begin{eqnarray*}
		 \pm c_{1}(X)   \pm c_{1}(Y) & \in &  
		 H^{2}(X, \ZZ)\oplus H^{2}(Y, \ZZ)  \\
	\pm c_{1}(X)   \pm c_{1}(Y)   \pm c_{1}(Z) 	 & \in &
	H^{2}(X, \ZZ)\oplus H^{2}(Y, \ZZ) 
			\oplus H^{2}(Z, \ZZ) 
	\end{eqnarray*}
	are monopole classes 
	on $X\# Y$ and $X\# Y\# Z$, respectively.
	Here the $\pm$ signs are arbitrary, and  independent
	of one another. 
		\end{cor}

\section{Obstructions to Einstein Metrics}\label{sec:nein}

When combined with the ideas of \cite{lric}, 
Corollary \ref{voila} immediately implies a non-existence
result for Einstein metrics. 

\begin{thm}\label{oui}
	Let $X$,  $Y$, and $Z$  be simply connected 
	symplectic $4$-manifolds
	with  $b_{+}\equiv 3 \bmod 4$. 
	Then $X\# Y\# k\overline{\CP}_{2}$ does not
	admit Einstein metrics if 
	$$k+4 \geq \frac{c_{1}^{2}(X) + c_{1}^{2}(Y)}{3}.$$
	Similarly, $X\# Y  \# Z \# k\overline{\CP}_{2}$ does not
	admit Einstein metrics if 
    $$k+8 \geq \frac{c_{1}^{2}(X) + c_{1}^{2}(Y)+c_{1}^{2}(Z)}{3}.$$
	\end{thm}
	
	\begin{proof}
		We begin by making the identification
		\begin{eqnarray*}
			 H^{2}(X \# k\overline{\CP}_{2},\ZZ)& = & 
			 H^{2}(X, \ZZ)\oplus \bigoplus_{j=1}^{k} H^{2}(\overline{\CP}_{2},
			 \ZZ) \\
			 & = & H^{2}(X, \ZZ)\oplus \ZZ^{\oplus k} , 
		\end{eqnarray*}
	in the process make a choice of generators $E_{j}$,
	$j=1,\ldots , k$ for the $k$ copies of $H^{2}(\overline{\CP}_{2},\ZZ)
	\cong \ZZ$. Now recall that there are self-diffeomorphisms
	of $X \# k\overline{\CP}_{2}$ which act trivially on 
		$H^{2}(X, \ZZ)$, but for which 
		$$E_{j}\mapsto \pm E_{j}$$
		for  any desired choice of signs. 
		Thinking of $X \# k\overline{\CP}_{2}$ as the 
		$k$-fold symplectic blow-up of $X$, and then 
		moving the symplectic structure via these 
		diffeomorphism, we thus obtain $2^{k}$ distinct
		symplectic structures on $X \# k\overline{\CP}_{2}$,
		with first Chern classes $$c_{1}= c_{1}(X)+\sum_{j=1}^{k}(\pm E_{j})$$
		for all possible choices of signs. 
		Applying Corollary \ref{voila} to 
		$(X \# k\overline{\CP}_{2})\# Y$ and 
		$(X \# k\overline{\CP}_{2})\# Y\# Z$,
		we thus conclude that all the classes of the form
		$\alpha+ \sum_{j=1}^{k} (\pm E_{j})$ are
		monopole classes, where  
		$$\alpha = 
		\left\{ 
		\begin{array}{ll}
			c_{1}(X) + c_{1}(Y), & \mbox{ if } M =  
			X\# Y\# k\overline{\CP}_{2},  \\
			c_{1}(X) + c_{1}(Y)+ c_{1}(Z), & \mbox{ if } M =  
			X\# Y\#  Z\# k\overline{\CP}_{2}. 
		\end{array}
		\right.$$

		Now, given any particular metric $g$ on $M$, let us
		make a new choice $\hat{E}_{j}=\pm E_{j}$ of generators 
		for our $k$ copies of $H^{2}(\overline{\CP}_{2},\ZZ)$
		in such a way that 
		$$\hat{E}_{j}\cdot \alpha^{+}\geq 0.$$
		The resulting monopole class
		$$a= \alpha + \sum_{j=1}^{k}\hat{E}_{j}$$
		then satisfies
		\begin{eqnarray*}
			(a^{+})^{2} & = & \left(\alpha^{+}+ \sum \hat{E}_{j}^{+}\right)^{2} 
			\\
			 & = &  (\alpha^{+})^{2}+ 2\sum \alpha^{+} \cdot
			  \hat{E}_{j}^{+}
			 +\left(\sum \hat{E}_{j}^{+} \right)^{2}  \\
			 & \geq  & (\alpha^{+})^{2}  \\
			 & \geq  & \alpha^{2}  
		\end{eqnarray*}
		Proposition \ref{est} therefore tells us that 
		any metric $g$ on $M$ satisfies
		\begin{equation}
			 	\frac{1}{4\pi^{2}}\int_{M}\left( \frac{s^{2}}{24} + 
2|W_{+}|^{2}\right) d\mu > \frac{2}{3} 
\alpha^{2} . 
			 	\label{eq:crux}
			 \end{equation}
(The  inequality is strict
because $a^{2}> (2\chi + 3\tau )(M)$, and this guarantees that    
$a$ is certainly  
not the first Chern class of a symplectic structure.) 

Now for any metric $g$ on our compact orientable $4$-manifold
$M$, we have the Gauss-Bonnet type formula \cite{bes,hit}
$$ (2\chi + 3\tau ) (M) = 	\frac{1}{4\pi^{2}}\int_{M}\left( 
\frac{s^{2}}{24} + 
2|W_{+}|^{2}-\frac{|\stackrel{\circ}{r}|^{2}}{2}\right) d\mu, $$
where $\stackrel{\circ}{r}=r-\frac{s}{4}g$  is the traceless 
Ricci tensor. If $g$ is Einstein, $\stackrel{\circ}{r}=0$,
and inequality (\ref{eq:crux}) then becomes
\begin{equation}
	(2\chi + 3\tau ) (M) >  \frac{2}{3} 
\alpha^{2}. 
	\label{eq:here}
\end{equation}

For  $M= X\# Y\# k\overline{\CP}_{2}$, we have 
$$\alpha^{2}= c_{1}^{2}(X) + c_{1}^{2}(Y),$$
and 
$$(2\chi + 3\tau ) (M)= c_{1}^{2}(X) + c_{1}^{2}(Y)-4-k.$$
Inequality (\ref{eq:here}) therefore asserts that 
a necessary condition for the existence of an Einstein metric
on $M$ is that 
$$c_{1}^{2}(X) + c_{1}^{2}(Y)-4-k > \frac{2}{3} [c_{1}^{2}(X) + 
c_{1}^{2}(Y)],
$$
or, in other words, that 
$$
\frac{c_{1}^{2}(X) + c_{1}^{2}(Y)}{3} > k+ 4. 
$$
By contraposition, this shows that there {\em cannot} be
an Einstein metric if
$$k+4 \geq \frac{c_{1}^{2}(X) + c_{1}^{2}(Y)}{3} ,$$
exactly as claimed. 

For $M=X\# Y\#  Z\# k\overline{\CP}_{2}$, we instead have
$$\alpha^{2}= 	c_{1}^{2}(X) + c_{1}^{2}(Y)+ c_{1}^{2}(Z),$$
and 
$$(2\chi + 3\tau ) (M)=	c_{1}^{2}(X) + c_{1}^{2}(Y)+ c_{1}^{2}(Z)
-8-k,$$
so that inequality (\ref{eq:here}) instead tells us that 
a necessary condition for 
 the existence of an Einstein metric
on $M$ is that 
$$c_{1}^{2}(X) + c_{1}^{2}(Y)+ c_{1}^{2}(Z)
-8-k > \frac{2}{3}[c_{1}^{2}(X) + c_{1}^{2}(Y)+ c_{1}^{2}(Z)],$$
or in other words that 
$$
\frac{c_{1}^{2}(X) + c_{1}^{2}(Y)+ c_{1}^{2}(Z)}{3} > k+ 8. 
$$
We thus conclude that there {\em cannot} be
an Einstein metric if
$$k+8 \geq \frac{c_{1}^{2}(X) + c_{1}^{2}(Y)+ c_{1}^{2}(Z)}{3} ,$$
 and this finishes the proof.  	
		\end{proof}
		
	In particular, setting $k=0$ yields: 
	
	\begin{cor}\label{non}
	Let $X$,  $Y$, and $Z$  be simply connected  symplectic $4$-manifolds
	with  $b_{+}\equiv 3 \bmod 4$. 
	If $c_{1}^{2}(X) + c_{1}^{2}(Y) \leq 12$, then 
	 $X\# Y$ does not
	admit Einstein metrics.
	Similarly, if $c_{1}^{2}(X) + c_{1}^{2}(Y)+c_{1}^{2}(Z) \leq 24$, 
	then 
	$X\# Y \#  Z$ does not
	admit Einstein metrics.
	\end{cor}

\section{Spin Examples}\label{sec:even}

In this section, we construct a sequence of $4$-dimensional spin 
manifolds which do not admit Einstein metrics, but nonetheless satisfy
(\ref{eq:strict}). Not only will we construct examples with an infinite
number of homeotypes,
but we will also show that 
each of these homeotypes carries an an infinite number of 
distinct smooth structures with the 
desired property. 

We begin with a collection of    
		symplectic spin manifolds constructed by Gompf \cite{gompf}.
		For arbitrary integers  $k  \geq 2$ and $m
		\geq 0$, Gompf showed that one can construct,
		by symplectic surgery methods, 
		a simply connected symplectic spin manifold
		with $(\chi, \tau ) = (24 k  + 4 m , -16 k )$. 
		In particular, setting $m = 2$, there
		is a simply connected symplectic spin $4$-manifold
		$X_{k }$
		with $(\chi, \tau ) = (24 k  + 8 , -16 k )$.
                   By the Minkowski-Hasse classification of quadratic forms, 
		this manifold thus has
		intersection form $-2k  {\bf E}_{8}\oplus (4k  +3){\bf H}$,
		and so is homeomorphic to $k (K3) \# (k +3)(S^{2}\times S^{2})$
		by the Freedman's theorem \cite{freedman}. 
		Notice that  $b_{+}(X_{k }) = 4k +3\equiv 3 \bmod 4$,
		and that  $c_{1}^{2}(X_{k })= (2\chi + 3\tau ) (X_{k })= 16$.

		It is perhaps worth remarking that most of these
		symplectic manifolds cannot be taken to be 
		complex surfaces, as most violate the Noether
		inequality $c_{1}^{2} \geq b_{+}-5$. On
		the other hand, we are free e.g. to take  
		$X_{4}$  to be a smooth complex 
		 hypersurface of tridegree $(4,4,2)$ in 
		 $\CP_{1}\times \CP_{1}\times \CP_{1}$.
		 In any case, it merely suffices for what follows that we 
		 choose one such symplectic manifold $X_{k}$
		 for each $k \geq 2$.

		 The other ingredient 
we will  need is a certain
 sequence of homotopy $K3$ surfaces.
Let $Y_{0}$ denote a Kummer surface, 
equipped with an elliptic
fibration $Y_{0}\to \CP_{1}$. Let
$Y_{\ell}$ be obtained from $Y_{0}$ by performing
a logarithmic transformation of order
$2\ell+1$ on one non-singular elliptic fiber of 
$Y_{0}$. The $Y_{\ell}$ are simply connected spin manifolds 
with $b_{+}=3$ and $b_{-}=19$, and hence
are homeomorphic to $K3$ by the 
Freedman classification. 
However, 
$Y_{\ell}$ is a K\"ahler surface
with $p_{g}=1\neq 0$, and so, for each $\ell$,
$\pm c_{1}(Y_{\ell})$
are   Seiberg-Witten basic class, with 
Seiberg-Witten invariant $\pm 1$. 
Moreover, $c_{1}(Y_{\ell})= 2\ell {\mathfrak f}$, where
$\mathfrak f$ is Poincar\'e dual to the multiple fiber
introduced by the logarithmic transformation \cite{bpv}. 

With these building blocks in hand, we now 
proceed to prove the following result:

\begin{thm} \label{spin} For any integer $n  \geq 4$, 
	the topological  spin manifold
	$$n (K3) \# (n+1) (S^{2}\times S^{2})$$
	admits infinitely many distinct smooth structures
	for which there cannot exist compatible Einstein metrics. 
	Moreover, this topological  manifold has 
	\begin{eqnarray*}
		\chi & = & 24n +4 , \\
		\tau  & = & -16n, 
	\end{eqnarray*}
	 and thus satisfies the strict
	Hitchin-Thorpe inequality $2\chi > 3 |\tau |$. 
	\end{thm}

	\begin{proof} Set $k=n-2$,  and 
		consider the simply connected $4$-manifolds 
		$$M_{k,\ell}=X_{k} \# Y_{0} \# Y_{\ell}.$$ Each of these 
		smooth oriented $4$-manifolds is homeomorphic to 
		$$(k +2) (K3) \# (k +3) (S^{2}\times S^{2}).$$ 
		However,
		since $b_{+}(X_{k})=4k + 3$ and $b_{+}(Y_{0})=b_{+}(Y_{\ell})=3$ 
		all reduce to $3 \bmod 4$, and since 
		$$c_{1}^{2}(X_{k})+c_{1}^{2}(Y_{0})+c_{1}^{2}(Y_{\ell})= 
		16+0+0 < 24,
		$$
		Corollary \ref{non} asserts that 
		none of these smooth manifolds $M_{k,\ell}$ 
		can admit an Einstein metric. 
		
		Now, for any fixed $k$,  observe that the sequence 
		$\{ M_{k,\ell} ~|~ \ell \in {\mathbb N}\}$
		must contain infinitely many
		different diffeotypes. Indeed, 
		Corollary \ref{voila} asserts that 
		$a=c_{1}(X_{k})+ c_{1}(Y_{\ell})$
		and $b=c_{1}(X_{k}) -c_{1}(Y_{\ell})$
		are both monopole classes on 
		$M_{k,\ell} =X_{k} \# Y_{0} \# Y_{\ell}$.  However, the  difference
		$a-b= 2c_{1}(Y_{\ell})$ 
		is divisible by $4\ell$, and  the bandwidth of $M_{k,\ell} $,
		as defined in  Definition \ref{band},  
		therefore satisfies 
		$$\BW (M_{k,\ell}  ) \geq 2\ell.$$
		Thus, for any fixed $k$,
		$$\sup_{\ell\in {\mathbb N}} ~\BW (M_{k,\ell})  =\infty,$$
		and it follows  that
		no individual $M_{k,\ell}$ has maximal bandwidth,
		 for any given $k$. 
		Hence
		$\{ M_{k,\ell }~|~ \ell \in {\mathbb N}\}$ 
		runs through 
		infinitely many different diffeotypes for each $k$,
		and the claim follows.
		\end{proof}

	\section{Non-Spin Examples}\label{sec:odd}

	We conclude this article by showing that  
Theorem \ref{oui} also provides a rich source 
of non-spin $4$-manifolds without Einstein metrics. 

We will begin with a sequence of  simply connected
symplectic manifolds constructed by Gompf 
\cite{gompf}.  Namely, for each integer $i \geq 2$, there is   
a compact simply connected symplectic $4$-manifold
$Z_{i }$ 
with  Todd genus $\frac{1}{4}(\chi + \tau )(Z_{i }) = i $
and $c_{1}^{2}(Z_{i }) = 8i -11$. 
The easiest thing to do 
with these examples is to apply the results of \cite{lric}, 
whereby one gets non-existence of Einstein metrics on
$Z_{i }\# k \overline{\CP}_{2}$ for $k\geq (8i  -11)/3$. Since the
resulting manifolds  are simply connected and non-spin, 
with $b_{+}=2i -1$ and $b_{-}=2i +10 +k$, Freedman's
theorem \cite{freedman} immediately yields 
	
\begin{prop}
	Let $(m,n)$ be a pair of natural numbers, 
	where $m$ is odd. If 
	$n > \frac{7}{3}m+8,$ there is
	a smooth structure on $m\CP_{2}\# n \overline{\CP}_{2}$
	for which there is no compatible Einstein metric. 
	\end{prop}
	
	Note that these examples satisfy the
	strict Hitchin-Thorpe inequality (\ref{eq:strict})
	provided that we also require that $n < 5m +4$;
	such examples therefore exist for any $m\neq 1$. 
 Moreover, many of these examples are actually
 homeomorphic to Einstein manifolds; cf. \cite{kot,lebweyl,lric}. 
 Of course,  this result does not pretend to be   
  optimal; e.g. sporadic improvements could certainly be made  by 
 exploiting the manifolds constructed by Stipsicz 
\cite{stip}.

	 	A much more interesting result is obtained, 
		however, by applying Theorem \ref{oui}
	to $Z_{2j}\# Y_{\ell}\# k \overline{\CP}_{2}$,
	where the $Y_{\ell}$ are the homotopy $K3$ surfaces 
	used in the previous
	section. If  $k > \frac{16}{3}j -8$, 
	these manifolds do not admit Einstein metrics. 
	On the other hand, for fixed $j,k$, the same
	reasoning used in the proof of Theorem \ref{spin}
	shows that 
	$$\lim_{\ell \to \infty}\BW (Z_{2j}\# Y_{\ell}\# k 
	\overline{\CP}_{2}) = \infty,$$
	so we obtain infinitely many distinct smooth structures
	for each  topological type. This proves the following:

	\begin{thm}\label{molti}
	Let $(m,n)$ be a pair of natural numbers
	with 
	 $m\equiv 2 \bmod 4$ and $m\geq 6$.  If 
	$n > \frac{7}{3}m+16,$ there are infinitely many 
	distinct  smooth structure on $m\CP_{2}\# n \overline{\CP}_{2}$
	 for which no  compatible Einstein metric exists. 
	\end{thm}
	
	Again, these manifolds will satisfy the strict
	Hitchin-Thorpe inequality if $n < 5m +4$,
	and such examples thus occur for each  
	allowed value of $m$.

	One can also apply the
	same method to $Z_{2j}\# R_{2,2}\# Y_{\ell}\# k 
	\overline{\CP}_{2}$, 
	where $R_{2,2}$ is the simply connected symplectic manifold 
	constructed by Gompf \cite{gompf}, 
	with  $b_{+}=3$ and $b_{-}=14$, and where 
	$k \geq \frac{16}{3}j -10$. By the same bandwidth argument, we obtain

\begin{thm}\label{lots}
	Let $(m,n)$ be a pair of natural numbers
	with 
	 $m\equiv 1 \bmod 4$ and $m\geq 9$.  If 
	$n > \frac{7}{3}m+12,$ there are infinitely many 
	distinct  smooth structure on $m\CP_{2}\# n \overline{\CP}_{2}$
	 for which no  compatible Einstein metric exists. 
	\end{thm}
	
	For each allowed value of $m$, we once again obtain 
	examples satisfying the  strict
	Hitchin-Thorpe inequality. But indeed, 
	many of these topological manifolds are known to 
	admit some {\em other} smooth structure for
	which a compatible Einstein metric {\em does}   exist. 
	
	\begin{xpl} For any  integer $p\geq 6$ with 
		$p\equiv 2 \bmod 4$,  set 	
	$m=p^{2}-3p+3$ and $n=3p^{2}-3p+1$. We then have
	$m\equiv 1\bmod 4$, $m \geq 9$,  and $n >  \frac{7}{3}m+12$;
	thus,   
	Theorem \ref{lots}
	therefore  asserts that  there are 
	infinitely many smooth structures on $m\CP_{2}\# n \overline{\CP}_{2}$
     for which no  Einstein metrics exists. However, this
	 homeotype is also realizable by the double 
	 branched  cover of $\CP_{2}$,
	 ramified over a smooth curve of degree
	 $2p$. This complex surface of general type contains no $(-1)$- or $(-2)$-curves, and so
	 has ample canonical line bundle
	 \cite{bpv}; the Aubin/Yau theorem
	 \cite{aubin,yau} thus predicts that it
	 admits an Einstein metric (which   can be constructed so as 
	  to also be K\"ahler).  
	 The situation is thus analogous to 
	 that of the $K3$ surface; for one smooth structure, 
	 there is an Einstein metric, but for infinitely many 
	 others, no such metric exists. 
	 \end{xpl}
		
	 Of course, the $K3$ surface has the remarkable additional 
	 property that
	 it admits  {\em only one} differentiable structure for
	 which there exists a compatible Einstein metric. 
	 In light of \cite{salvetti}, however, no such uniqueness statement 
	 holds for many  of the topological manifolds under 
	 discussion; cf. \cite{cat,kot}. Nonetheless, one might 
	 speculate that  
	 Einstein metrics can   exist  only for a 
	  {\em finite number} 
	 of smooth structures on any given topological $4$-manifold.
	 Perhaps this difficult question will become a fruitful
	 topic for future
	 research. 

\bigskip
\noindent
{\bf Acknowledgment.} 
We would like to  express our deep gratitude to Professors Stefan Bauer and  Mikio Furuta
for  helping us  come to grips with  the key features of  the  Bauer-Furuta invariant.

  \end{document}